\documentclass[titlepage,twoside,12pt]{article}
\usepackage{amssymb}
\usepackage{amsfonts}
\begin{document}

{\LARGE \bf The Necessary Structure of \\ \\ Congruences in Free Semigroups} \\ \\


{\bf Elem\'{e}r E ~Rosinger} \\ \\
{\small \it Department of Mathematics \\ and Applied Mathematics} \\
{\small \it University of Pretoria} \\
{\small \it Pretoria} \\
{\small \it 0002 South Africa} \\
{\small \it eerosinger@hotmail.com} \\ \\

{\bf Abstract} \\

A characterization of congruences in free semigroups is presented. \\ \\

{\bf 1. Semigroups} \\

In this, and the next four sections, we recall for convenience of notation and terminology several basic concepts and
constructions related to semigroups, [1,2]. The characterization of congruences in free semigroups is presented in section
6. This seemingly new result has recently had important applications which are published elsewhere. \\

A {\it semigroup} is a structure $( S, \star )$, where $S$ is a nonvoid set and $\star : S \times S \longrightarrow S$ is
a binary operation on $S$ which is {\it associative}, that is, it satisfies the condition \\

(1.1)~~~ $ \star ( u, \star ( v, w ) ) = \star ( \star ( u, v ), w ),~~~ u, v, w \in S $ \\

Here we recall that it is customary to denote $\star ( u, v )$ simply by $u \star v$, for $u, v \in S$. Consequently, the
above associativity condition is equivalent with \\

(1.2)~~~ $ u \star ( v \star w ) = ( u \star v ) \star w,~~~ u, v, w \in S $ \\

We note that semigroups need not always have to be commutative, or have neutral elements. \\

A fundamental concept which can relate semigroups to one another is recalled now. Given two semigroups $( S, \star )$ and
$( T, \diamond )$, a mapping $f : S \longrightarrow T$ is called a {\it homomorphism}, if and only if \\

(1.3)~~~ $ f ( u \star v ) = f ( u ) \diamond f ( v ),~~~ u, v \in S $ \\

If such a homomorphism $f$ is injective, then it is called a {\it monomorphism}. In case a homomorphism $f$ is surjective,
then it is called an {\it epimorphism}. And if a homomorphism $f$ is surjective and it also has an inverse mapping $f^{-1}
: T \longrightarrow S$ which is again a homomorphism, then it is called an {\it isomorphism}. It follows that isomorphic
semigroups are identical for all purposes, and they only differ in the notation of their elements, or of the semigroups
themselves. \\ \\

{\bf 2. Free Semigroups} \\

Free semigroups are of fundamental importance since, as indicated in Proposition 2.1. below, all semigroups can be
obtained from them in a natural manner. \\

Given any nonvoid set $E$, we denote by $E^+$ the set of all finite sequences \\

(2.1)~~~ $ a_1, a_2, \ldots , a_n $ \\

where $n \geq 1$ and $a_i \in E$, with $1 \leq i \leq n$. Further, we define on $E^+$ the binary operation $\ast$ as
follows \\

(2.2)~~~ $ ( a_1, a_2, \ldots , a_n ) \ast ( b_1, b_2, \ldots , b_m ) =
                                  a_1, a_2, \ldots , a_n, b_1, b_2, \ldots , b_m $ \\

that is, simply by the juxtaposition of sequences in  $E^+$. Then it is easy to see that \\

(2.3)~~~ $ ( E^+, \ast ) $ is a semigroup \\

since $\ast$ is obviously associative. However, if $E$ has at least two elements, then $\ast$ is clearly not commutative.
Also, the semigroup $( E^+, \ast )$ does not have a neutral element, regardless of the number of elements in $E$. \\

A consequence of the associativity of $\ast$ is that the elements of the semigroup $E^+$ can be written in the following
form which is alternative to (2.1), namely \\

(2.4)~~~ $ a_1 \ast a_2 \ast \ldots \ast a_n $ \\

where $n \geq 1$ and $a_i \in E$, with $1 \leq i \leq n$. \\

The semigroup $( E^+, \ast )$ is called the {\it free semigroup on} $E$, and the meaning of that term will result from the
two important properties presented next. \\

First, the mapping \\

(2.5)~~~ $ E \ni a \stackrel{i}\longmapsto a \in E^+ $ \\

is {\it injective}, and obviously, it is never surjective. \\

Second, the free semigroup $( E^+, \ast )$ on $E$ has the following {\it universality property}. Given any semigroup $( S,
\star )$ and any mapping $j : E \longrightarrow S$, there exists a {\it unique} homomorphism $f : E^+ \longrightarrow S$
such that \\

(2.6)~~~ $ j = f \circ i $ \\

or in other words, the diagram commutes

\begin{math}
\setlength{\unitlength}{1cm}
\thicklines
\begin{picture}(13,6)

\put(0,2.4){$(2.7)$}
\put(1.1,5){$E$}
\put(2,5.1){\vector(1,0){7}}
\put(9.5,5){$E^+$}
\put(5.5,5.4){$i$}
\put(1.7,4.5){\vector(1,-1){3.3}}
\put(9.5,4.5){\vector(-1,-1){3.3}}
\put(5.4,0.6){$S$}
\put(2.8,2.5){$j$}
\put(8.1,2.5){$\exists !~~f$}

\end{picture}
\end{math} \\

The above universality property of semigroups has an important immediate consequence, namely \\

{\bf Proposition 2.1.} \\

Every semigroup is the homomorphic image of a free semigroup. \\

{\bf Proof.} \\

Let $( S, \star )$ be a semigroup, then we can take $E = S$ and $j = id_S$ in (2.7), and obtain the commutative diagram

\begin{math}
\setlength{\unitlength}{1cm}
\thicklines
\begin{picture}(13,7)

\put(1.1,5){$S$}
\put(2,5.1){\vector(1,0){7}}
\put(9.5,5){$S^+$}
\put(5.5,5.4){$i$}
\put(1.7,4.5){\vector(1,-1){3.3}}
\put(9.5,4.5){\vector(-1,-1){3.3}}
\put(0,2.5){$(2.8)$}
\put(5.4,0.6){$S$}
\put(2.8,2.5){$id_S$}
\put(8.1,2.5){$\exists !~~f$}

\end{picture}
\end{math} \\

However, in view of (2.6), obviously $f$ is surjective. Thus $S$ is indeed the homomorphic image of the free semigroup
$S^+$. \\ \\

{\bf 3. Quotient Constructions} \\

Let $E$ be a nonvoid set and $\approx$ an equivalence relation on $E$. Then the {\it quotient set} \\

(3.1)~~~ $ E / \approx $ \\

is defined as having the elements given by the {\it cosets} \\

(3.2)~~~ $ ( a )_\approx = \{~ b \in E ~~|~~ b \approx a ~\},~~~ a \in E $ \\

thus each coset $( a )_\approx$ is the set of all elements $b \in E$ which are equivalent with $a$ with respect to
$\approx$. The coset $( a )_\approx$ is also called the {\it equivalence class} of $a$ with respect to the equivalence
relation $\approx$. It follows that the mapping \\

(3.3)~~~ $ i_\approx : E \ni a \longmapsto ( a )_\approx \in E / \approx $ \\

is surjective, and it is called the {\it canonical quotient mapping}. \\

A useful way to obtain equivalence relations on any given set $E$ is through the construction called {\it transitive
closure}. Namely, given any family $( \equiv_i )_{i \in I}$ of symmetric binary relations on $E$, then we
define the equivalence relation $\equiv$ on $E$ as follows. If $a, b \in E$, then \\

(3.4)~~~ $ a \equiv b $ \\

holds, if and only if $a = b$, or there exist $c_0, c_1, c_2, \ldots , c_n \in E,~ i_1, i_2, i_3, \ldots , i_n \in I$, with $n \geq 1$, such that \\

(3.5)~~~ $ a ~=~ c_0 ~\equiv_{i_1}~ c_1 ~\equiv_{i_2}~ c_2 ~\equiv_{i_3}~ \ldots ~\equiv_{i_n}~ c_n ~=~ b $ \\

It is convenient to identify any binary relation $r$ on a nonvoid set $E$ with the subset of $E \times E$ given by \\

(3.6)~~~ $ r = ~\{~ ( a, b ) \in E \times E ~~|~~ a ~r~ b ~\} $ \\

It follows easily that the transitive closure of the family $( \equiv_i )_{i \in I}$ of symmetric binary relations on $E$ is the same with the
transitive closure of the corresponding single symmetric binary relation $\equiv $ on $E$ given by \\

(3.7)~~~ $ \equiv ~=~ \cup_{i \in I} \equiv_i $ \\

In general, for a symmetric binary relation $r$ on $E$, we shall denote by $r^{tc}$ its transitive  closure. \\

An alternative and equivalent way to construct quotient spaces is through {\it partitions}. Given a partition of $E$ by
the family ${\cal E} = ( E_i )_{i \in I}$ of subsets of $E$. Then one can associate with it an equivalence relation
$\approx_{\cal E}$ on $E$, defined for $a, b \in E$, by \\

(3.8)~~~ $ a \approx_{\cal E} b ~~~\Longleftrightarrow~~~ \exists~~ i \in I ~:~ a, b \in E_i $ \\

Obviously, in this case we have for $a \in E$ and $i \in I$ \\

(3.9)~~~ $ a \in E_i ~~~\Longleftrightarrow~~~ E_i = ( a )_{\approx_{\cal E}} $ \\

in other words, the equivalence class $( a )_{\approx_{\cal E}}$ of $a$ with respect to $\approx_{\cal E}$ is precisely
the set $E_i$ in the partition ${\cal E}$ to which $a$ belongs. Consequently \\

(3.10)~~~ $ E / \approx_{\cal E} ~=~ \{~ E_i ~~|~~ i \in I ~\} $ \\ \\

{\bf 4. Congruences} \\

Let $( S, \star )$ be any semigroup. An equivalence relation $\approx$ on $S$ is called a {\it congruence} on $( S,
\star )$, if and only if it is compatible with the semigroup operation $\star$ in the following sense \\

(4.1)~~~ $ u \approx v ~~~\Longrightarrow~~ u \star w \approx v \star w,~~ w \star u \approx w \star v $ \\

for all $u, v, w \in S$. \\

The importance of such a congruence is that the resulting quotient $S / \approx $ of $S$ leads again to a {\it semigroup},
namely \\

(4.2)~~~ $ ( S, \star ) / \approx ~~~=~~~ (\, S / \approx,\, \diamond ) $ \\

where the binary operation $\diamond$ on $S / \approx$ is defined by \\

(4.3)~~~ $ ( u )_\approx \diamond ( v )_\approx = ( u \star v )_\approx,~~~ u, v \in S $ \\

also, the canonical quotient mapping, see (3.4) \\

(4.4)~~~ $ S \ni u \longmapsto ( u )_\approx \in S / \approx $ \\

is a {\it surjective homomorphism}, thus an {\it epimorphism}. \\

Furthermore, let $( S, \star )$ and $( T, \diamond )$ be two semigroups and $f : S \longrightarrow T$ a morphism between
them. Then the binary relation on $S$ given by \\

(4.5)~~~ $ ker f = \{~ ( u, v ) \in S \times S ~~|~~ f ( u ) = f ( v ) ~\} $ \\

is a congruence on $( S, \star )$, and there exists a monomorphism $g : ( S, \star ) / ker f \longrightarrow ( T,
\diamond )$, such that, see (3.3) \\

(4.6)~~~ $ f = g \circ i_{ker f} $ \\

which means that the diagram commutes \\

\begin{math}
\setlength{\unitlength}{1cm}
\thicklines
\begin{picture}(13,6)

\put(0,2.5){$(4.7)$}
\put(1.1,5){$S$}
\put(2,5.1){\vector(1,0){7.3}}
\put(9.9,5){$T$}
\put(5.5,5.4){$f$}
\put(1.7,4.5){\vector(1,-1){3.3}}
\put(6.5,1.3){\vector(1,1){3.3}}
\put(5.1,0.6){$S / ker f$}
\put(2.5,2.5){$i_{ker f}$}
\put(8.1,2.5){$g$}

\end{picture}
\end{math} \\

Given now a partition ${\cal S} = ( S_i )_{i \in I}$ of $S$, then in view of (3.8) - (3.10), it leads to an equivalence relation
$\approx_{\cal S}$ on $S$. Indeed, in view of (4.1), it is obvious that the equivalence relation $\approx_{\cal S}$
on $S$ will be a {\it congruence} on $( S, \star )$, if and only if, for every $i \in I,~ u, v \in S_i,~ w \in S$, we
have \\

(4.8)~~~ $ ( u \star w )_{\approx_{\cal S}} = ( v \star w )_{\approx_{\cal S}},~~~
               ( w \star u )_{\approx_{\cal S}} = ( w \star v )_{\approx_{\cal S}} $ \\ \\

{\bf 5. The Case of Commutative Semigroups} \\

Let $( E^+, \ast )$ be the free semigroup generated by the nonvoid set $E$, see (2.3). Given two sequences, see (2.4) \\

(5.1)~~~ $ a_1 \ast \ldots \ast a_n,~~~ b_1 \ast \ldots \ast b_m \in E^+ $ \\

we define \\

(5.2)~~~ $ a_1 \ast \ldots \ast a_n ~\approx~ b_1 \ast \ldots \ast b_m $ \\

if and only if the two sequences are the same, or differ by a permutation of their elements. Then obviously $\approx$ is
an equivalence relation on $E^+$ which is also a congruence on $( E^+, \ast )$. It follows that \\

(5.3)~~~ $ ( E^+, \ast ) / \approx $ \\

is a commutative semigroup. \\

The commutative version of Proposition 2.1. is \\

{\bf Proposition 5.1.} \\

If $( S, \ast )$ is a commutative semigroup, then it is the homomorphic image of $( S^+, \ast ) / \approx$~ through the
mapping \\

(5.4)~~~ $ S^+ / \approx ~\ni ( a_1 \ast \ldots \ast a_n )_\approx \longmapsto a_1 \ast \ldots \ast a_n \in S $ \\

where $a_1, \ldots , a_n \in S$. \\ \\

{\bf 6. Congruences on Free Semigroups} \\

Given a nonvoid set $E$, let $( E^+, \ast )$ be the corresponding free semigroup, as in (2.3) above. Further, let
$\approx$ be an equivalence relation on $E^+$. Then it follows easily that $\approx$ is a congruence on $E^+$, if and only
if, for every $a \in E,~ u, v \in E^+$, we have \\

(6.1)~~~ $ u \approx v ~~~\Longrightarrow~~~ a u ~\approx~ a v,~~~ u a ~\approx~ v a $ \\

 We associate with $\approx$ the binary relation $\approx_{min}$ on $E^+$ as follows. Given $u, v \in E^+$, then $u
 \approx_{min} v$, if and only if $u \approx v$, while neither of the following two relations holds \\

(6.2)~~~ $ u = a u\,',~ v = a v\,',~ u\,' ~\approx~ v\,' $ \\

for certain $a \in E,~ u\,', v\,' \in E^+$, or \\

(6.3)~~~ $ u = u\,'' b,~ v = v\,'' b ,~ u\,'' ~\approx~ v\,'' $ \\

for certain $b \in E,~ u\,'', v\,'' \in E^+$. \\

Obviously \\

(6.4)~~~ $ ~\approx_{min}~ \subseteq ~\approx~ $ ~and~ $ ~\approx_{min}~ $ is symmetric \\

also, we have the following weak version of reflexivity for $\approx_{min}$ \\

(6.5)~~~ $ a ~\approx_{min}~ a,~~~ a \in E $ \\

Moreover, if $u, v \in E^+$ and $u \approx v$, then in view of (6.2), (6.3), we have \\

(6.6)~~~ $ u = a_1 \ldots a_n u\,' b_1 \ldots b_m,~~~ v = a_1 \ldots a_n v\,' b_1 \ldots b_m,~~~
                                                            u\,' \approx_{min} v\,' $ \\

for suitable $a_1 \ldots a_n, b_1 \ldots b_m \in E$, with $n, m \geq 0$, and $u\,', v\,' \in E^+$. \\

In view of (6.6), we introduce two concepts. A given arbitrary binary relation $r \subseteq E^+ \times E^+$ is called {\it
invariant}, if an only if for every $a \in E,~ u, v \in E^+$, we have, see (6.1) \\

(6.10)~~~ $ u ~r~ v ~~~\Longrightarrow~~~ a u ~r~ a v,~~~ u a ~r~ v a $ \\

Further, the binary relation $r^{ic} \subseteq E^+ \times E^+$, called {\it invariant closure} of $r$, is defined as
follows. For $u, v \in E^+$ we have $u ~r^{ic}~ v$, if and only if \\

(6.11)~~~ $ u = a_1 \ldots a_n u\,' b_1 \ldots b_m,~~~ v = a_1 \ldots a_n v\,' b_1 \ldots b_m,~~~
                                                                                        u\,' ~r~ v\,' $ \\

for suitable $a_1 \ldots a_n, b_1 \ldots b_m \in E$, with $n, m \geq 0$, and $u\,', v\,' \in E^+$. \\

{\bf Lemma 6.1.} \\

For any equivalence relation $~\approx~$ on $E^+$ we have \\

(6.12)~~~ $ ~\approx_{min}~ \subseteq ~\approx^{tc}_{min}~ \subseteq ~\approx~ \subseteq ~\approx^{ic}_{min}~ $ \\

{\bf Proof.} \\

The inclusion $\approx ~\subseteq~ \approx^{ic}_{min}$ follows clearly from (6.6). \\

Given now $u, v \in E^+$, with $u ~\approx^{tc}_{min}~ v$, then (3.4), (3.5) lead to $u = v$, or alternatively to \\

(6.13)~~~ $ u = z_1 ~\approx_{min}~ \ldots ~\approx_{min}~ z_n = v $ \\

for some $z_1, \ldots , z_n \in E^+$. And then (6.4) implies \\

(6.14)~~~ $ u = z_1 ~\approx~ \ldots ~\approx~ z_n = v $ \\

thus $u ~\approx~ v$.

\hfill $\Box$ \\

And now, we obtain the characterization in \\

{\bf Theorem 6.1.} \\

The equivalence relation $\approx$ is a congruence on the free semigroup $( E^+, \ast )$, if and only if it is the
invariant closure of $\approx_{min}$, namely \\

(6.15)~~~ $ \approx ~~=~~ \approx^{ic}_{min} $ \\

in which case we have \\

(6.16)~~~ $ ~\approx_{min}~ \subseteq ~\approx^{tc}_{min}~ \subseteq ~\approx~ ~~=~~ \approx^{ic}_{min} $ \\

{\bf Proof.} \\

If the equality (6.15) holds, then $\approx$ is obviously a congruence, in view of (6.1), (6.11). \\

Conversely, let us assume that the equivalence relation $\approx$ is a congruence on the free semigroup $( E^+, \ast )$. The (6.1), (6.11) and (6.12)
imply (6.15). \\


\begin{thebibliography}{99}

\bibitem{} Bourbaki N : Elements of Mathematics, Algebra I, \\ Chapters 1-3. Springer-Verlag, New York, September
1998, \\ ISBN-13: 9783540642435

\bibitem{} Howie J M : Fundamentals of Semigroup Theory. Calderon, Oxford, 2003

\end{thebibliography}
\end{document}